\documentclass 
{amsart}
\usepackage{amsmath}
\usepackage{amsthm}
\usepackage{amscd}
\usepackage{amsxtra}
\usepackage{amssymb}
\usepackage[all]{xy}
\usepackage{boxedminipage}
\usepackage{ulem}
\usepackage{stmaryrd}
\usepackage[dvipdfm]{graphicx,color}

\newtheorem{theorem}{Theorem}[section]
\newtheorem{lemma}[theorem]{Lemma}

\newtheorem{corollary}[theorem]{Corollary}

\newtheorem*{DKK31}{{\color{red}{Corrected Version}} of Proposition 3.1. of \cite{DKK}}
\newtheorem*{DKK32}{{\color{red}{Corrected Version}} of Proposition 3.2. of \cite{DKK}}
\newtheorem*{DKKT}{{\color{red}{Corrected Version}} of the identity of \cite{DKK}}

\newtheorem*{H-S}{Theorem (Hovey-Sadofsky)}

\theoremstyle{definition}
\newtheorem{definition}[theorem]{Definition}
\newtheorem{example}[theorem]{Example}

\newtheorem{def-pr}[theorem]{Definition-Proposition}
\newtheorem{def-th}[theorem]{Definition-Theorem}
\newtheorem{def-cor}[theorem]{Definition-Corollary}

\newcommand{\IM}{ \operatorname{Im} }   

\newcommand{\Hom}{ \operatorname{Hom} }

\newcommand{\Epi}{ \operatorname{Epi} }
\newcommand{\Mono}{ \operatorname{Mono} }

\newcommand{\lcm}{ \operatorname{lcm} }

\newcommand{\ord}{ \operatorname{ord} }

\newcommand{\disc}{ \operatorname{disc} }

\newcommand{\Z}{{ \mathbb Z}}

\newcommand{\F}{{ \mathbb F}}

\newcommand{\N}{{ \mathbb N}}

\newcommand{\C}{{ \mathbb C}}

\newcommand{\Ocal}{ {\mathcal{O}} }

\newcommand{\Ab}{ {\mathfrak{Ab}} }

\newcommand{\dotbox}{\hbox to 1em{\hss.\hss}}

\begin{document}

\title{On the random variable $\N \ni l \mapsto \gcd(l,n_1) \gcd(l, n_2) \cdots \gcd(l, n_k) \in \N$}
\author{Norihiko Minami}
\address{Omohi College, Nagoya Institute of Technology, Gokiso, Showa-ku, Nagoya 466-8555}
\email{nori@nitech.ac.jp}

\begin{abstract}
For natural numbers $n_1, \ldots, n_k$, evaluating the moments of the random variable
\begin{equation*} 
\begin{split}
X: \Omega := \left\{ 1, 2, \ldots, \lcm(n_1,n_2,\ldots,n_k)  \right\} &\to \N   \\
   l &\mapsto \gcd(l,n_1) \gcd(l, n_2) \cdots \gcd(l, n_k)
\end{split}
\end{equation*}
by some purely elementary method, we obtain a series of identitities of elementary
number theory.
The special case of average was originally considered by Deitmar-Koyama-Kurokawa \cite{DKK}, who studied this case
by an analytic consideration of some zeta function.

This average turns out to be an extremely important quantity in the Soul\'e type zeta functions of $\F_1$-schemes,
and we show this average is nothing but the invariant $\mu(A)$ 
of an abelian group $A := \prod_{j=1}^k ( \Z/ n_j \Z)$, defined by
\begin{equation*}
\mu (A) = \sum_{a\in A} \frac{1}{ | a | }.
\end{equation*}
\end{abstract}
\maketitle

\section{Introduction and main results}

In his talk at JAMI 2009 on March 24, Kurokawa presented a rather mysterious looking identidity
of elementary number theory:
\begin{equation} \label{FI}
\frac{1}{n}\sum_{k=1}^n \gcd(n,k) = \prod_{p|n} \left(  1 + \left( 1 - \frac{1}{p} \right) \ord_p(n) \right)
\end{equation}

\begin{example}  $n=12 = 2^2\cdot 3$:
\begin{multline*}
\frac{1}{12} ( 1 + 2 + 3 + 4 + 1 + 6 + 1 + 4 +  3 + 2 + 1 + 12)  = \frac{40}{12} =  \frac{10}{3} \\
=  2\cdot  \frac{5}{3} = \left(  1 + \left( 1 - \frac{1}{2} \right)\cdot 2 \right)
\cdot \left(  1 + \left( 1 - \frac{1}{3} \right)\cdot 1 \right)
\end{multline*}
\end{example}

In fact, this is a special case (the case $k=1$) of identities considered by \cite{DKK}.

\begin{DKKT}
For $n_1, n_2, \ldots, n_k \in \N$,
\begin{multline} \label{DKKI}
\frac{1}{{\color{red}{\lcm(n_1,n_2,\ldots,n_k)}}}\sum_{l =1}^{{\color{red}{\lcm(n_1,n_2,\ldots,n_k)}}} \gcd(l,n_1) \gcd(l,n_2)\cdots \gcd(l, n_k)  \\
= 
\prod_{ p \big| \lcm(n_1,n_2,\ldots,n_k)}  
\left[
{\color{red}{p^{\nu_{p,0} + \cdots + \nu_{p,k-1} }}} + \left( 1 - \frac{1}{p} \right) 
\sum_{j=0}^{k-1} p^{\nu_{p,{\color{red}{0}}} + \cdots + \nu_{p,j}} \sum_{\mu=\nu_{p,j}}^{\nu_{p,j+1}-1} 
p^{{\color{red}{(k-j)\nu_{p,j}-\mu}}}
\right]
\end{multline} 
Here, for each prime $p \big| \lcm(n_1,n_2,\ldots,n_k)$, 
\begin{gather*}
\{ \nu_{p,1}, \nu_{p,2}, \ldots, \nu_{p,k-1} , \nu_{p,k} \}
= \{ \ord_p(n_1),  \ord_p(n_2), \ldots,   \ord_p(n_{k-1}),  \ord_p(n_k)  \}   \\
\nu_{p,0} := 0 \leq \nu_{p,1} \leq \nu_{p,2} \leq \ldots \leq \nu_{p,k-1} \leq \nu_{p,k}
\end{gather*}
\end{DKKT}

\begin{example}  $n_1 = 6, n_2 = 4,\ k=2$:
\begin{multline*}
\frac{1}{\lcm(6,4)} (  1\cdot 1 + 2\cdot 2 + 3\cdot 1 +  2\cdot 4  + 1\cdot 1 + 6\cdot 2
+ 1\cdot 1 + 2\cdot 4 + 3\cdot 1 + 2\cdot 2 + 1\cdot 1 + 6\cdot 4 )
  \\
= \frac{1}{12}( 1 + 4 + 3 + 8 + 1 + 12 + 1 + 8 + 3 + 4 + 1 + 24) = \frac{70}{12}=\frac{35}{6} 
\\
= \frac{7}{2}\cdot \frac{5}{3} 
= \left( 2 + \frac{1}{2} ( 1\times 1 + 2 \times 1 ) \right)
\cdot \left( 1 +  \frac{2}{3}( 1\cdot 0 + 1\cdot 1 ) \right)
 \\
= \left(  2^1 + \left( 1 - \frac{1}{2} \right) 
\left(  2^0 \sum_{\mu=0}^{1-1} 2^{(2-0)\cdot 0 - \mu}
+ 2^1 \sum_{\mu=1}^{2-1} 2^{(2-1)\cdot 1 - \mu} \right) \right)   \\
\times 
\left(  3^0 + \left( 1 - \frac{1}{3} \right) 
\left(  3^0 \sum_{\mu=0}^{0-1} 3^{(2-0)\cdot 0 - \mu}
+ 3^0 \sum_{\mu=0}^{1-1} 3^{(2-1)\cdot 0 - \mu} \right) \right) 
\end{multline*}
\end{example}

\cite{DKK} tried to obtain \eqref{DKKI} by studying some zeta function of Igusa type,
and Kurokawa said he is not aware of any elementary proof even for \eqref{FI}.

Now the purpose of this paper is to give a purely elementary proof of generalizations of
\eqref{DKKI} from the view point of elementary probability theory.
Fixing a finite abelian group $A := \prod_{j=1}^k (\Z/n_j\Z)$ with 
$n_1, n_2, \ldots, n_k \in\N$, we would like to understand the random variable:

\begin{equation} \label{old}
\begin{split}
\tilde{X}(A): \tilde{\Omega} := \N &\to \N   \\
       l &\mapsto  \# \Big| \Hom_{\Ab}( A, \Z/l\Z ) \Big|    
= \gcd(l,n_1) \gcd(l, n_2) \cdots \gcd(l, n_k) ,
\end{split}
\end{equation}
where $\Ab$ is the category of abelian groups.

Although $\tilde{\Omega} = \N$ is an infinite set, we would like to regard it being equipped with
the \lq\lq homogeneous measure\rq\rq.  For this purpose, we observe:
\begin{multline*}
l \equiv l'\ \mod \lcm(n_1,n_2,\ldots,n_k)  \\
\implies   
\gcd(l,n_1) \gcd(l, n_2) \cdots \gcd(l, n_k)
= \gcd(l',n_1) \gcd(l', n_2) \cdots \gcd(l', n_k)
\end{multline*}
Thus, instead of \eqref{old}, we may equally consider the following random variable:
\begin{equation} \label{new}
\begin{split}
X(A) : \Omega := \left\{ 1, 2, \ldots, \lcm(n_1,n_2,\ldots,n_k)  \right\} &\to \N   \\
   l &\mapsto \gcd(l,n_1) \gcd(l, n_2) \cdots \gcd(l, n_k) ,     
\end{split}
\end{equation}
where $\Omega= \left\{ 1, 2, \ldots, \lcm(n_1,n_2,\ldots,n_k)  \right\}$ is equipped with the homogeneous measure.
Then the identity \eqref{DKKI} is nothing but a convenient formula to evaluate the average $E[X(A)]$ of
the random variable \eqref{new}.   However, from the view point of elementary probability theory, it is very natural to
seek for similar convenient formulae for the variance $V[X(A)] = E[X(A)^2] - E[X(A)]^2$ and even \lq\lq higher\rq\rq invariants.  

Our Main Theorem 
offers such formulae for the continuous version $E[X(A)^w]\ (w\in \C)$. Their
special cases $w\in\N$ are nothing but the moments of the randowm variable $X(A)$, and the simplest case $w=1$
is nothing but the identity \eqref{DKKI}.

\begin{theorem} \label{MT}
For a finite abelian group $A = \prod_{j=1}^k (\Z/n_j\Z)$ with $n_1, n_2, \ldots, n_k \in \N,$ and $w\in \C$,
{
\begin{multline} \label{MTI}
\hspace{-4mm}E[X(A)^w] =
\frac{1}{
\lcm(n_1,n_2,\ldots,n_k)
}\sum_{l =1}^{
{\lcm(n_1,n_2,\ldots,n_k)
}} \left( \gcd(l,n_1) \gcd(l,n_2)\cdots \gcd(l, n_k) \right)^w \\
\hspace{-20mm}= 
\prod_{ p \big| \lcm(n_1,n_2,\ldots,n_k)}  
\left[
p^{ ( \nu_{p,0} + \cdots + \nu_{p,k-1} )w + \nu(p,k)(w-1) } 
\right.
\\
\left.\hspace{50mm} 
+ \left( 1 - \frac{1}{p} \right) 
\sum_{j=0}^{k-1} p^{ ( \nu_{p,
0} + \cdots + \nu_{p,j} )w } \sum_{\mu=\nu_{p,j}}^{\nu_{p,j+1}-1} p^{
(k-j)\nu_{p,j}w-\mu}
\right]
\\
\hspace{-20mm}= 
\prod_{ p \big| \lcm(n_1,n_2,\ldots,n_k)}  
\left[
p^{ ( \nu_{p,0} + \cdots + \nu_{p,k-1} )w + \nu(p,k)(w-1) } 
\right.
\\
\left.\hspace{50mm} 
 + \left( 1 - \frac{1}{p} \right) 
\sum_{j=0}^{k-1} p^{ \left( \nu_{p,
0} + \cdots + \nu_{p,j-1} + (k-j+1)\nu_{p,j} \right)w } \sum_{\mu=\nu_{p,j}}^{\nu_{p,j+1}-1} p^{
-\mu}
\right]
\end{multline} 
}
Here, for each prime $p \big| \lcm(n_1,n_2,\ldots,n_k)$, 
\begin{gather*}
\{ \nu_{p,1}, \nu_{p,2}, \ldots, \nu_{p,k-1} , \nu_{p,k} \}
= \{ \ord_p(n_1),  \ord_p(n_2), \ldots,   \ord_p(n_{k-1}),  \ord_p(n_k)  \}   \\
\nu_{p,0} := 0 \leq \nu_{p,1} \leq \nu_{p,2} \leq \ldots \leq \nu_{p,k-1} \leq \nu_{p,k}
\end{gather*}
\end{theorem}


Our purely elementary proof of Theorem~\ref{MT} turns out to be a matter of triviality, when restricted to
the original identity \eqref{FI}.
For our general case, when $w\in \N$, we have the following:

\begin{lemma} \label{winN}
For $w\in \N$ and a finite abelian group $A = \prod_{j=1}^k\left(\Z/n_j\Z\right)$,
consider the $w$ fold product group $A^w = \prod_{j=1}^k\left(\Z/n_j\Z\right)^w$. Then we have the following:
\begin{equation}
E[ X(A)^w ] = E[ X(A^w) ]
\end{equation}
\end{lemma}

Thus, the esense of Theorem~\ref{MT}, at least when $w\in\N$, would be an effective evaluation
of 
$E[ X(A) ]$ for a general finite abelian group $A$. Now, the following theorem does the job for us:

\begin{theorem} \label{mu-evaluation}
For any finite abelian group $A$,
\begin{equation} \label{muA}
E[ X(A) ]
= \sum_{a \in A} \frac{1}{ |a| } ,
\end{equation}
where $|a|$ stands for the order of an element $a\in A$.
\end{theorem}

Motivated by Theorem~\ref{mu-evaluation}, we make the following definition:

\begin{definition}
For a finite abelian group $A$, define
\begin{equation}
\mu(A) := \sum_{a \in A} \frac{1}{ |a| } 
\end{equation}
\end{definition}

Thus, Theorem~\ref{mu-evaluation} may be restated as:
\begin{equation} \label{restated}
E[ X(A) ] = \mu(A)
\end{equation}

Applying Theorem~\ref{mu-evaluation} (or \eqref{restated}), 
we would like to obtain, when $w\in\N$, more user-frinedly evaluation of $E[ X(A)^w ]$ 
than Theorem~\ref{MT}.
For this purpose, applying the Chinese Remaindar Theorem, 
we express a finite abelian group 
\begin{equation} \label{CRT}
A = \prod_{j=1}^k\left(\Z/n_j\Z\right) 
= \prod_p  \prod_{j=1}^k \left(\Z/p^{\ord_p(n_j)}\Z\right),
\end{equation}
in the following form:
\begin{equation} \label{convenient}
A = \prod_{j=1}^k\left(\Z/n_j\Z\right) 
= \prod_p \prod_{i=1}^{h_p} ( \Z/ p^i )^{m_{p,i}} ,
\end{equation}
where 
\begin{gather*} 
h_p := \max \left\{ \ord_p(n_j) \mid j = 1, \cdots, k \right\}, \\
m_{p,i} := \# \{ j = 1, \cdots, k \mid  \ord_p(n_j) = i \} \quad ( 1\leq i \leq h_p),
\end{gather*}

In \eqref{convenient}, we immediately see the following for $1\leq l\leq h_p$:
\begin{equation} \label{change}
\begin{split}
&\quad \# \left\{ g \in \prod_{i=1}^{h_p} ( \Z/ p^i )^{m_{p,i}} \Big| p^l g = 0 \right\}    \\ 
&= p^{\sum_{i=1}^{l} im_{p,i}  + l \sum_{i=l+1}^{h_p} m_{p,i}  }  
= p^{\sum_{i=1}^{l-1} im_{p,i}  + l \sum_{i=l}^{h_p} m_{p,i}  }  \\
&\quad \# \left\{ g \in \prod_{i=1}^{h_p} ( \Z/ p^i )^{m_{p,i}} \Big| |g| = p^l \right\}  \\ 
&=  \# \left\{ g \in \prod_{i=1}^{h_p} ( \Z/ p^i )^{m_{p,i}} \Big| p^l g = 0 \right\} 
-  \# \left\{ g \in \prod_{i=1}^{h_p} ( \Z/ p^i )^{m_{p,i}} \Big| p^{l-1} g = 0 \right\}  \\
&= p^{\sum_{i=1}^{l-1} im_{p,i}  + l \sum_{i=l}^{h_p} m_{p,i}  } -
p^{\sum_{i=1}^{l-1} im_{p,i}  + (l-1) \sum_{i=l}^{h_p} m_{p,i}  }   \\
&= \left( p^{ \sum_{i=l}^{h_p} m_{p,i} } - 1 \right) 
p^{( \sum_{i=1}^{l-1} im_{p,i}  + (l-1) \sum_{i=l}^{h_p} m_{p,i} ) } 
\end{split}
\end{equation}

Now, Theorem~\ref{mu-evaluation} gives us the following user-frinedly version of 
Therem~\ref{MT} for $w\in\N$:

\begin{theorem} \label{user-friendly}
For $w\in\N$, $A = \prod_{j=1}^k\left(\Z/n_j\Z\right) 
= \prod_p \prod_{i=1}^{h_p} ( \Z/ p^i )^{m_{p,i}}$,
as in \eqref{convenient},
\begin{multline} \label{friendly}
\begin{split}
E[ X( A )^w ] &= \prod_p \mu \left(  \prod_{i=1}^{h_p} ( \Z/ p^i )^{wm_{p,i}} \right)  \\  
&= \prod_p \left(  1 + \sum_{l=1}^{h_p}
 \left( p^{ w\sum_{i=l}^{h_p} m_{p,i} } - 1 \right) 
p^{w\left(\sum_{i=1}^{l-1} im_{p,i}  + (l-1) \sum_{i=l}^{h_p} m_{p,i} \right) } \frac{1}{p^l} \right) \\
&= \prod_p \left(  1 + \sum_{l=1}^{h_p}
 \left( p^{ w\sum_{i=l}^{h_p} m_{p,i} } - 1 \right) 
p^{\left(w\sum_{i=1}^{l-1} im_{p,i}  + w(l-1)\sum_{i=l}^{h_p} m_{p,i} - l) \right)} \right)
\end{split}
\end{multline}
\end{theorem}

In Section 2, we shall prove Theorem~\ref{MT}.  In Section 3, we shall prove 
Theorem~\ref{mu-evaluation} and Theorem~\ref{user-friendly}.

As far as the original identity \eqref{FI} concerns, there is yet another kind of transparent proof due to 
Alain Connes, using the Euler $\phi$ function. 
We have contained this proof of Connes as Appendix 1. 
We would like to express our gratidute to Professor Connes for allowing us to contain his argument in this paper.

Since the original Igusa zeta calculations of \eqref{DKKI} are not correctly stated in \cite{DKK},  we presented 
corrected computations in 
Appendex 2 for readers' convenience.

We note another kind of generalizations of \eqref{FI} is given in \cite{KO}, which
is once again proved using some zeta functions of Igusa type.  Even for this and some generalizations, we can
offer an elementary proof \cite{M1}.  

Motivated by results of this paper and \cite{M1}, we defined and studied some multivarible {\it deformed} 
zeta function of $\F_1$-scheme of Hurewicz-Igusa-type \cite{M2},
which generalizes both the zeta functions studied in \cite{DKK} \cite{KO} and the log derivative
of the modified Soul\'e type zeta function \cite{CC}.

Amongst of all, the following Theorem~\ref{zeta-formula} is proved in \cite{M2}, which manifestly shows
the importance of $E[ X(A) ] = \mu(A)$, 
studied in this paper and \cite{DKK}, in the study of zeta functions of $F_1$-scheme.


\begin{theorem} \label{zeta-formula}
For a Noetherian $F_1$-scheme $X$, let $\zeta_X(s)$ be the (generalized) Soul\'e zeta function \cite{S} \cite{CC}, 
and let $\zeta^{\disc}_X(s)$ be the modified zeta function, both of which were defined and studued in \cite{CC}.
Then, there are some entire functions $h_1(s), h_2(s)$ s.t.
\begin{equation*} \label{The expression}
\begin{split}
&\quad \zeta_X(s) = e^{h_1(s)} \zeta^{\disc}_X(s)   \\
&\hspace{-20mm}= e^{h_2(s) } \prod_{x\in X}
\left(
\left(  \prod_{j=0}^{n(x)} (s-j)^{ \left( - \binom{ n(x) }{j}   (-1)^{ n(x) - j } \right) } \right)^{   
{\color{blue}\mu\left( \prod_j \Z/ m_j(x)\Z \right) } }
\right) ,
\end{split}
\end{equation*}
where, for each $x\in X$, 
$
\Ocal_x^{\times} 
= \Z^{n(x)}\times {\color{blue}\prod_j \Z/ m_j(x)\Z}.$ 
\end{theorem}

The basic ideas of the results in this paper were obtained during the author's stay at JAMI2009, Johns Hopkins University, in March 2009,
and presented at NCGOA2009, Vanderbilt University, in May 2009, and at the Fall Meeting of the Mathematical Society of Japan
at Osaka University, in September 2009.

The author would like to express his gratitude to 
Jack Morava, Takashi Ono, Steve Wilson, and Guoliang Yu for their hospitalities.
The author also would like to express his gratitude to 
Katia Concani, Alain Connes, and Nobushige Kurokawa for their work and encouragement.

\section{Proof of Theorem~\ref{MT}}

Whereas \cite{DKK} used some multivariable Igusa type zeta funciions for 
groups,
we use some 
ring structure along the line of \cite{M1}:
For a finite 
ring with $n$ elements $R = \{ l_i \in R \mid 1\leq i\leq n\}$ and its ideals $I_1, I_2, \ldots, I_k$, set
\begin{equation*} 
Z_{(R; I_1, I_2, \ldots, I_k)}
(w) := \frac{1}{ | R | } \sum_{  l \in R  } \left( \big| R/ (l + I_1) \big| \cdot \big| R/ (l + I_2) \big| \cdots
\big| R/ (l + I_k) \big|  \right)^w
\quad 
(w\in \C)
\end{equation*}
We wish to understand this, because
\begin{multline*}
Z_{ \left( \Z/ \lcm(n_1,n_2,\ldots,n_k) \Z; (n_1), (n_2), \ldots, (n_k) \right)   }
(w) 
\\
=
\frac{1}{
\lcm(n_1,n_2,\ldots,n_k)
}\sum_{l =1}^{
{\lcm(n_1,n_2,\ldots,n_k)
}} \left( \gcd(l,n_1) \gcd(l,n_2)\cdots \gcd(l, n_k) \right)^w 
\quad 
(w\in \C)
\end{multline*}

Of course, we have an elementary probability theoretical interpretation:
For the random variable
\begin{align*}
X: \Omega := R &\to \N   \\
      l &\mapsto  
\big| R/ (l + I_1) \big| \cdot \big| R/ (l + I_2) \big| \cdots
\big| R/ (l + I_k) \big|,
\end{align*}
where $\Omega = R$ is equipped with the homogeneous measure,
\begin{equation*}
E[X^w] = \frac{1}{ | R | } \sum_{  l \in R  } \left( \big| R/ (l + I_1) \big| \cdot \big| R/ (l + I_2) \big| \cdots
\big| R/ (l + I_k) \big|  \right)^w
= 
Z_{(R; I_1, I_2, \ldots, I_k)}
(w) 
\end{equation*}
%


Then the following easy observation turns out to play an essential role in our proof of Theorem~\ref{MT}:

\begin{lemma} \label{product}
If there is a finite ring decomposition $R = \prod_{i} R_i$ and the corresponding decompositions
of the ideals
$I_j = \prod_i I_{i,j}$ with $I_{i,j} \subset R_i$ for $1\leq j\leq k$. Then
\begin{equation*} \label{decomposition}
Z_{(R; I_1, I_2, \ldots, I_k)}(w)
= \prod_i Z_{(R_i; I_{i,1}, I_{i,2}, \ldots, I_{i,k})}(w)
\hspace{30mm}
\qed
\end{equation*}
\end{lemma}

\begin{corollary}
For $n_1, n_2, \ldots, n_k \in \N$ and each prime $p \big| \lcm(n_1,n_2,\ldots,n_k)$, set:
\begin{gather*}
\{ \nu_{p,1}, \nu_{p,2}, \ldots, \nu_{p,k-1} , \nu_{p,k} \}
= \{ \ord_p(n_1),  \ord_p(n_2), \ldots,   \ord_p(n_{k-1}),  \ord_p(n_k)  \}   \\
\nu_{p,0} := 0 \leq \nu_{p,1} \leq \nu_{p,2} \leq \ldots \leq \nu_{p,k-1} \leq \nu_{p,k}  .
\end{gather*}
Then,
\begin{multline}  \label{Euler-product}
\frac{1}{
\lcm(n_1,n_2,\ldots,n_k)
}\sum_{l =1}^{
{\lcm(n_1,n_2,\ldots,n_k)
}} \left( \gcd(l,n_1) \gcd(l,n_2)\cdots \gcd(l, n_k) \right)^w 
\\
=
\prod_{ p \big| \lcm(n_1,n_2,\ldots,n_k)}  
\left[
\frac{1}{p^{ \nu_{p,k} }} \sum_{l=1}^{ p^{ \nu_{p,k} } }
\left( \gcd(l,p^{ \nu_{p,1} }) \gcd(l,p^{ \nu_{p,2} })\cdots \gcd(l, p^{ \nu_{p,k} }) \right)^w \right]
\end{multline}
\end{corollary}

\begin{proof}
This follows immediately from Lemma~\ref{product} and the Chinese Remainder Theorem.
\end{proof}

Now the corresponding $p$-factor random variable
\begin{equation}
\begin{split}
X_p: [1, p^{ \nu_{p,k} } ] := \{ l\in \N \mid 1\leq l\leq   p^{ \nu_{p,k} }   \}  &\to \N   \\
l &\mapsto 
\gcd(l,p^{ \nu_{p,1} }) \gcd(l,p^{ \nu_{p,2} })\cdots \gcd(l, p^{ \nu_{p,k} }) 
\end{split}
\end{equation}
factors as
\begin{equation*}
X_p
: [1, p^{ \nu_{p,k} } ]
\xrightarrow{ \ord_p }  [0,  \nu_{p,k} ]
\xrightarrow{ Y_p } \N .
\end{equation*}
Here,
\begin{gather*}
Y_p:  [0, \nu_{p,k} ]  
= \{ \nu_{p,k}  \} \coprod
\{ \mu_j\in \N \mid 0\leq j< k,\ \nu_{p,j}\leq \mu_j < \nu_{p, j+1} \}
\to \N   \\
 \nu_{p,k} \mapsto    p^{\nu_{p,1}}\cdots
 p^{\nu_{p,k-1}  }
  p^{\nu_{p,k} }
=p^{\nu_{p,0} + \nu_{p,1} + \cdots + \nu_{p,k-1} + \nu_{p,k} }  \\
 \hspace{26mm} \mu_j \mapsto p^{\nu_{p,0}}p^{\nu_{p,1}}\cdots
 p^{\nu_{p,j-1}  }  p^{\nu_{p,j} }\cdot \left(  p^{\nu_{p,j} } \right)^{k-j}
= p^{\nu_{p,0} + \nu_{p,1} + \cdots + \nu_{p,j-1} + \nu_{p,j} } 
 p^{(k-j)\nu_{p,j} } ,
\end{gather*}
and observe, for $\mu\in [0, \nu_{p,k}]$,
\begin{equation}
\big| \ord_p^{-1}(\mu) \big| =
\begin{cases}
1\quad &\text{if}\ \mu= \nu_{p,k}  \\
p^{\nu_{p,k} - \mu} - p^{\nu_{p,k} - \mu -1}\ &\text{if}\ \mu<\nu_{p,k}
\end{cases}
\end{equation}

Therefore,
\begin{equation} \label{p-factor}
\begin{split}
&\quad
\frac{1}{p^{ \nu_{p,k} }} \sum_{l=1}^{ p^{ \nu_{p,k} } }
\left( \gcd(l,p^{ \nu_{p,1} }) \gcd(l,p^{ \nu_{p,2} })\cdots \gcd(l, p^{ \nu_{p,k} }) \right)^w
\\
&= E[X_p^w]
= \frac{1}{ p^{ \nu_{p,k} } } \sum_{l=1}^{ p^{ \nu_{p,k} } } X_p(l)^w  
= \frac{1}{ p^{ \nu_{p,k} } } \sum_{\mu\in [0, \nu_{p,k} ] } \big| \ord_p^{-1}(\mu) \big| Y_p(\mu)^w
\\
&= \frac{1}{ p^{ \nu_{p,k} } } 
\left[
1\cdot \left( 
p^{\nu_{p,0} + \nu_{p,1} + \cdots + \nu_{p,k-1} + \nu_{p,k} }  
 \right)^w
 \right.
 \\
&\hspace{20mm}
\left.
+ \sum_{j=0}^{k-1} \sum_{\mu_j = \nu_{p,j}}^{\nu_{p,j+1}-1}
\left( p^{\nu_{p,k} - \mu_j} - p^{\nu_{p,k} - \mu_j -1} \right) \cdot
\left(
p^{\nu_{p,0} + \nu_{p,1} + \cdots + \nu_{p,j-1} + \nu_{p,j} } 
 p^{(k-j)\nu_{p,j} } 
\right)^w
\right]
\\
&= 
p^{ ( \nu_{p,0} + \cdots + \nu_{p,k-1} )w + \nu(p,k)(w-1) } 
 + \left( 1 - \frac{1}{p} \right) 
\sum_{j=0}^{k-1} p^{ ( \nu_{p,
0} + \cdots + \nu_{p,j} )w } \sum_{\mu_j=\nu_{p,j}}^{\nu_{p,j+1}-1} p^{
(k-j)\nu_{p,j}w-\mu_j}
\\
&= 
p^{ ( \nu_{p,0} + \cdots + \nu_{p,k-1} )w + \nu(p,k)(w-1) } 
 + \left( 1 - \frac{1}{p} \right) 
\sum_{j=0}^{k-1} p^{ \left( \nu_{p,
0} + \cdots + \nu_{p,j-1} + (k-j+1)\nu_{p,j} \right)w } \sum_{\mu_j=\nu_{p,j}}^{\nu_{p,j+1}-1} p^{
-\mu_j}
\end{split}
\end{equation}

Combining \eqref{Euler-product} and \eqref{p-factor},
Theorem~\ref{MT} has been proven.
\qed

\section{Proof of Theorem~\ref{mu-evaluation}, and Theorem~\ref{user-friendly}}




\begin{proof}[Proof of Theorem~\ref{mu-evaluation}]
We work in the category $\Ab$ of abelian groups, and for abelian groups $A$ and $B$,
we denote by  $\Hom_{\Ab}(A,B),\ \Epi_{\Ab}(A, B),\ \Mono_{\Ab}(A, B),$ the set of
all homomorphisms, all epimorhisms, and all monomorphisms, from $A$ to $B$, respectively.
Then, for a finite abelian group $A = \prod_{j=1}^k\left(\Z/n_j\Z\right)$  and its
cyclic subgroup $C \subset A$, we note
\begin{equation} \label{mono}
\Big|  \Mono_{\Ab}( C, \Z/l\Z) \Big| =
\begin{cases}
\phi ( | C | ) \quad &\text{if}\ | C | \mid l \\
0 \quad &\text{if}\ | C | \nmid l
\end{cases}
\end{equation}
where $\phi$ is Euler's totient function. 
Furthermore, 
\begin{equation} \label{cyclic}
\Big| \left\{ h \in \Hom_{\Ab}(\Z, A ) \mid h(\Z) = C \right\} \Big| = \phi ( | C | )
\end{equation}
Then the formula \eqref{muA} is obtained in the following order:
\begin{align*}
&\qquad E[ X(A) ] \\
&:=
\frac{1}{{\color{red}{\lcm(n_1,n_2,\ldots,n_k)}}}\sum_{l =1}^{{\color{red}{\lcm(n_1,n_2,\ldots,n_k)}}} \gcd(l,n_1) \gcd(l,n_2)\cdots \gcd(l, n_k)  \\
&=
\frac{1}{ n_1\cdot n_2\cdots n_k}\sum_{l =1}^{ n_1\cdot n_2\cdots n_k} 
\gcd(l,n_1) \gcd(l,n_2)\cdots \gcd(l, n_k)  \\
&= \frac{1}{ |A| }\sum_{l=1}^{ |A| }  \Big|  \Hom_{\Ab}( A, \Z/l\Z ) \Big|
= \frac{1}{ |A| }\sum_{l=1}^{ |A| }  \Big|  \Hom_{\Ab}( \Z/l\Z, A ) \Big|   
\\
&= \frac{1}{ |A| }\sum_{l=1}^{ |A| } \sum_{cyclic\ C \subset A}
\Big|  \Epi_{\Ab}( \Z/l\Z, C) \Big|   
= \frac{1}{ |A| } \sum_{cyclic\ C \subset A} \sum_{l=1}^{ |A| }
\Big|  \Epi_{\Ab}( \Z/l\Z, C) \Big|  \\
&= \frac{1}{ |A| } \sum_{cyclic\ C \subset A} \sum_{l=1}^{ |A| }
\Big|  \Mono_{\Ab}( C, \Z/l\Z) \Big|   
\overset{\because \eqref{mono}}=  \frac{1}{ |A| } \sum_{cyclic\ C \subset A} \frac{ |A|}{|C|}\phi(|C|)  \\
&= \sum_{cyclic\ C \subset A} \frac{ \phi ( |C| ) }{ |C| }  
\overset{\because \eqref{cyclic}}= \sum_{h \in \Hom(\Z, A )} \frac{1}{ | h(1) | }  = \sum_{a\in A} \frac{1}{ |a| }
\end{align*}
\end{proof}

\begin{proof}[Proof of Theorem~\ref{user-friendly}]
When $A$ is as in \eqref{convenient},
\begin{equation} \label{convenient-w}
A^w = \prod_{j=1}^k\left(\Z/n_j\Z\right)^w 
= \prod_p \prod_{i=1}^{h_p} ( \Z/ p^i )^{wm_{p,i}}
\end{equation}
By Lemma~\ref{winN} and Theorem~\ref{mu-evaluation}, it is very easy to see
\begin{equation} \label{EulerFactorization}
E[X(A)^w] = E[X(A^w)] = \mu (A^w ) = \prod_p \mu \left(  \prod_{i=1}^{h_p} ( \Z/ p^i )^{wm_{p,i}} \right) .
\end{equation}
Thus, it suffices to evaluate $\mu \left(  \prod_{i=1}^{h_p} ( \Z/ p^i )^{wm_{p,i}} \right)$.
For this purpose, we consider the version of \eqref{change} with $A^w$
for $1\leq l\leq h_p$:
\begin{equation} \label{change-w}
\begin{split}
&\quad \# \left\{ g \in \prod_{i=1}^{h_p} ( \Z/ p^i )^{wm_{p,i}} \Big| p^l g = 0 \right\} \\
&= p^{\sum_{i=1}^{l} iwm_{p,i}  + l \sum_{i=l+1}^{h_p} wm_{p,i}  }  
= p^{\sum_{i=1}^{l-1} iwm_{p,i}  + l \sum_{i=l}^{h_p} wm_{p,i}  }  \\ 
&= p^{w\sum_{i=1}^{l-1} im_{p,i}  + wl \sum_{i=l}^{h_p} m_{p,i}  }  \\
&\quad \# \left\{ g \in \prod_{i=1}^{h_p} ( \Z/ p^i )^{wm_{p,i}} \Big| |g| = p^l \right\}  \\
&= p^{\sum_{i=1}^{l-1} iwm_{p,i}  + l \sum_{i=l}^{h_p} wm_{p,i}  } -
p^{\sum_{i=1}^{l-1} iwm_{p,i}  + (l-1) \sum_{i=l}^{h_p} wm_{p,i}  }   \\
&= \left( p^{ \sum_{i=l}^{h_p} wm_{p,i} } - 1 \right) 
p^{( \sum_{i=1}^{l-1} iwm_{p,i}  + (l-1) \sum_{i=l}^{h_p} wm_{p,i} ) }  \\
&= \left( p^{ w\sum_{i=l}^{h_p} m_{p,i} } - 1 \right) 
p^{( w\sum_{i=1}^{l-1} im_{p,i}  + w(l-1) \sum_{i=l}^{h_p} m_{p,i} ) } 
\end{split}
\end{equation}
Now the claim follows immediately by applying Theorem~\ref{mu-evaluation} and \eqref{change-w}.
\end{proof}

\section{Appendix 1: Connes' proof of \eqref{FI}}

Conne's proof of \eqref{FI} makes full use of the Euler totient function
$\phi$, which allows us to derive the Euler product decomposition in \eqref{FI}
without resorting to the Chinese Remainder Theorem.  The key observation 
to relate the Euler totient function $\phi$ is the following fact, which was used 
in the proof of Lemma 5.7 in \cite{CC}:

\begin{lemma} \label{Connes} For $n, m\in\N$,
\begin{equation}
\gcd(n,m) = \sum_{{d | n,\ d | m}} \phi(d)
\end{equation}
\end{lemma}

\begin{proof}
Observe that
$
\gcd(n,m) = \Big| \Hom_{\text{Abelian group}}(\Z/n, \Z/m) \Big|  .
$
Then, for each homomorphism $h\in \Hom_{\text{Abelian group}}(\Z/n, \Z/m)$, we can associate
$d\in\N$ such that
\begin{equation} \label{image}
\IM h \cong \Z/d\quad (d | n,\ d| m)
\end{equation}
When $d$ is fixed, there are exactly $\phi(d)$ homomorphisms $h$ satifying \eqref{image}.
So the claim follows.
\end{proof}

\begin{proof}[Connes' Proof of \eqref{FI}]
\begin{align*}
&\qquad \frac{1}{n}\sum_{k=1}^n \gcd(n,k) = \frac{1}{n}\sum_{k=1}^n \sum_{{d | n,\ d | k}} \phi(d)
= \frac{1}{n} \sum_{d | n} \Big| \{ k\in \N \mid 1\leq k\leq n,\ d | k \} \Big| \phi(d)
\\
&= \frac{1}{n}\sum_{d | n} \frac{n}{d} \phi(d) = \sum_{d | n}  \frac{ \phi(d) }{d} 
= \sum_{ ( j_p )_{ p | n} \in \prod_{ p \mid n} [0, \ord_p(n) ] } \frac{ \phi\left( \prod_{p | n} p^{j_p} \right) }{  \prod_{p | n} p^{j_p} }
\\
&= \sum_{ ( j_p )_{ p | n} \in \prod_{ p \mid n} [0, \ord_p(n) ] }  \prod_{p | n} \frac{ \phi( p^{j_p} ) }{ p^{j_p} }
= \prod_{ p \mid n} \sum_{ j_p =0}^{ \ord_p(n)  } \frac{ \phi( p^{j_p} ) }{ p^{j_p} }
\\
&= \prod_{ p \mid n} \left( 1 + \sum_{ j_p =1}^{ \ord_p(n)  } \frac{ \phi( p^{j_p} ) }{ p^{j_p} } \right)
= \prod_{p|n} \left(  1 + \left( 1 - \frac{1}{p} \right) \ord_p(n) \right) 
\end{align*}
\end{proof}

Using Lemma~5.8  of \cite{CC}, Connes' proof may be generalized to prove \eqref{DKKI}.

\section{Appendix 2: The corrected Igusa-type zeta computation of \cite{DKK}}

For a finitely generated abeian group $A$, \cite{DKK} defined and studied its absolute zeta function of Igusa type $\zeta^I(s,A)$:
\begin{equation*}
\zeta^I(s,A) := \sum_{m=1}^{+\infty} \frac{ \big| \Hom_{\Ab}(A, \Z/m\Z) \big| }{m^s}
\end{equation*}
\cite{DKK} tried to consider the general case
$A = \Z^r \times (\Z/n_1\Z) \times \cdots \times (\Z/n_k\Z)$
by computing this in two different ways, corresponding to:
\begin{itemize}
\item $m = \prod_{p\mid m} p^{\ord_p(m)}$
\item $m = \nu\cdot \lcm(n_1,n_2,\ldots,n_k) + l,\ \nu\geq 0, 1\leq l\leq  \lcm(n_1,n_2,\ldots,n_k)$
\end{itemize}

Unfortunately, the computations of \cite{DKK} contain some mistake.
Such is the case, we offer the {\color{red}{corrected}} computations for readers' convenience:

\begin{DKK31}
\begin{align*}
&\quad Z^{\text{group}}\left( s;  \Z^r \times (\Z/n_1\Z ) \times \cdots \times ( \Z/n_k\Z ) \right)  \\
&= \zeta(s-r) \prod_{p\mid n}  
\left[
p^{{\color{red}{ (\nu_{p,1} + \cdots + \nu_{p,k-1} + \nu_{p,k} )+ (r-s)\nu_{p,k} }}} 
\right.\\
&\left.\hspace{30mm}
+ ( 1 - p^{r-s} ) 
\sum_{j=0}^{k-1} p^{\nu_{p,{\color{red}{0}}} + \cdots + \nu_{p,j}} \sum_{\mu_j=\nu_{p,j}}^{\nu_{p,j+1}-1} 
p^{{\color{red}{(k-j)\nu_{p,j}+(r-s)\mu_j}}}
\right]
\end{align*}
Here, for each prime $p \mid n$, 
\begin{gather*}
\{ \nu_{p,1}, \nu_{p,2}, \ldots, \nu_{p,k-1} , \nu_{p,k} \}
= \{ \ord_p(n_1),  \ord_p(n_2), \ldots,   \ord_p(n_{k-1}),  \ord_p(n_k)  \}   \\
\nu_{p,0} := 0 \leq \nu_{p,1} \leq \nu_{p,2} \leq \ldots \leq \nu_{p,k-1} \leq \nu_{p,k}
\end{gather*}
\end{DKK31}

\begin{DKK32}
\begin{align*}
&\quad \zeta^I
\left( s;  \Z^r \times (\Z/n_1\Z ) \times \cdots \times ( \Z/n_k\Z ) \right)  \\
&= {\color{red}{\lcm(n_1,n_2,\ldots,n_k)}}^{r-s} \\
&\hspace{3mm}\times \sum_{l =1}^{{\color{red}{\lcm(n_1,n_2,\ldots,n_k)}}} \gcd(l,n_1) \gcd(l,n_2)\cdots \gcd(l, n_k) 
\times \zeta \left( s-r, \frac{l}{ {\color{red}{\lcm(n_1,n_2,\ldots,n_k)}} } \right),
\end{align*}
where $\zeta(s, q ) := \sum_{m=0}^{\infty} (m+q)^{-s}\ (\Re(q)>0)$ is the {\it Hurwitz zeta function},
which has a simple pole at $s=1$  with residue $1$.
\end{DKK32}

Both expressions imply $\zeta^I
\left( s;  \Z^r \times (\Z/n_1\Z ) \times \cdots \times ( \Z/n_k\Z ) \right)$ 
has a meromorphic continuation to $\C^r$
Then, as in \cite{DKK}, noticing that 
\begin{equation*}
\zeta^{I}
\left( s;  \Z^r \times (\Z/n_1\Z ) \times \cdots \times ( \Z/n_k\Z ) \right) 
\end{equation*}
has a simple pole at $s=r+1,$ and evaluating the residue at $s=r+1$ in two different ways,
we obtain the identity \eqref{DKKI}. 


\end{document}